\documentclass{amsart}
\usepackage{amsfonts}

\newtheorem{theorem}{Theorem}[section]

\theoremstyle{definition}

\theoremstyle{remark}

\numberwithin{equation}{section}
\theoremstyle{plain}

\copyrightinfo{2001}{enter name of copyright holder}
\input{tcilatex}

\begin{document}
\title[Almost Quarter-Pinched Manifolds]{Classification of Almost
Quarter-Pinched Manifolds}
\author{Peter Petersen}
\address{520 Portola Plaza\\
Dept of Math UCLA\\
Los Angeles, CA 90095}
\email{petersen@math.ucla.edu}
\urladdr{http://www.math.ucla.edu/\symbol{126}petersen}
\thanks{}
\author{Terence Tao}
\email{tao@math.ucla.edu}
\urladdr{http://www.math.ucla.edu/\symbol{126}tao}
\date{}
\subjclass{}
\keywords{}

\begin{abstract}
We show that if a simply connected manifold is almost quarter pinched then
it is diffeomorphic to a CROSS or sphere.
\end{abstract}

\maketitle

\section{Introduction}

The goal of this paper is to use the new results of Brendle and Schoen in 
\cite{Brendle-SchoenI}, \cite{Brendle-SchoenII} to sharpen some older
results about almost quarter-pinched manifolds studied by Berger and
Abresch-Meyer in \cite{Berger}, \cite{Abresh-MeyerII}. The main goal is to
show that there can't be exotic spheres which are almost quarter pinched.
More precisely, we prove

\begin{theorem}
\label{main} There exist $\varepsilon \left( n\right) >0$ so that any simply
connected $n$-dimensional Riemannian manifold $\left( M,g\right) $ with
sectional curvatures in $\left[ 1/4-\varepsilon ,1\right] $ is diffeomorphic
to a sphere or compact rank one symmetric space (CROSS).
\end{theorem}

In \cite{Berger}, Berger established this result in even dimensions, but
with the weaker conclusion that $M$ was either \emph{homeomorphic} to a
sphere or diffeomorphic to a CROSS (see also \cite{Durumeric} for other
results in this direction.) To prove his result, Berger considered a
sequence of almost quarter-pinched manifolds, then using Cheeger-Gromov
compactness theory he passed to a limit that was weakly quarter-pinched in
the sense of comparison theory. It does not seem that a similar strategy
would work for Theorem \ref{main}; the difficulty is that there is no
obvious way to use comparison curvatures to get information about isotropic
curvatures on manifolds that only have weakly defined curvature tensors.
Instead we apply the Ricci flow to each of the manifolds in the sequence and
show that in the limit we get a Ricci flow where we have nice metrics for
positive time. These are not necessarily quarter-pinched (as this property
is not preserved by Ricci flow), but they do have nonnegative isotropic
curvatures when we add a factor of $\mathbb{R}^{2}$ to the metric. Thus we
can use the classification in \cite{Brendle-SchoenII} to understand what the
candidates for limit manifolds are. Finally it is worth mentioning that
there is at least one example of an exotic sphere with positive sectional
curvature (see \cite{Petersen-Wilhelm}.)

All results in Riemannian geometry which are not specifically referenced can
be found in \cite{Petersen} and similarly for results on the Ricci Flow in 
\cite{Chow-Knopf} and \cite{Chow-Lu-Ni}.

The first author would like to thank Frederick Wilhelm for bringing the
above question to our attention. The second author is supported by NSF grant
DMS-0649473 and a grant from the Macarthur Foundation.

\section{Proof of Theorem \protect\ref{main}}

Assume by contradiction that Theorem \ref{main} failed in some fixed
dimension $n$. Then one could find a sequence of simply connected $n$%
-dimensional Riemannian manifolds $\left( M_{i},g_{i}\right) $ such that $%
1/4-1/i\leq \mathrm{sec}\leq 1$, and such that none of the $M_{i}$ were
diffeomorphic to a sphere or a CROSS.

In even dimensions it is a classical result of Klingenberg that such
manifolds have injectivity radius $\geq \pi .$ In odd dimensions this is far
more subtle and was settled by Abresch and Meyer in \cite{Abresch-MeyerI}.
The fact that there is a uniform lower bound for the sectional curvature
shows that the diameter is also bounded, by Myers' theorem. Finally, the
sectional curvature bounds imply that the curvature tensors of $g_{i}$ are
uniformly bounded. By the standard local existence theory for Ricci flow (%
\cite{HamiltonII}, \cite{HamiltonIII}, and \cite{Chow-Knopf}), we can thus
run the Ricci flow for a fixed amount of time $\left[ 0,t_{0}\right] $ for
each of these metrics $(M_{i},g_{i})$, with the curvature tensor uniformly
bounded in $i$ on this time interval.

Putting all this together, we obtain a family of Ricci flows $\left(
M_{i},g_{i}\left( t\right) ,t\in \left[ 0,t_{0}\right] \right) $ with
uniformly bounded geometry. Hamilton's extension of the Cheeger-Gromov
compactness theorem then guarantees us that some subsequence converges to a
Ricci flow $\left( M,g\left( t\right) ,t\in \left[ 0,t_{0}\right] \right) $ (%
\cite{HamiltonIV} and \cite{Chow-Lu-Ni}.); the convergence is only in the $%
C^{1,\alpha },$ $\alpha <1$ sense at time $t=0$, but is in the $C^{\infty }$
sense for $0<t\leq t_{0}$.

We cannot expect $g\left( t\right) $ to be quarter-pinched for $t>0,$ and
when $t=0$ the metric isn't sufficiently smooth that this makes sense.
However, for small $t$ the metrics $g_{i}\left( t\right) $ have sectional
curvature $\geq 1/4-1/i-C\left( n\right) t$ for an absolute constant $C>0$
(see \cite[Proposition 2.5]{Rong}). As the metrics $g_{i}\left( t\right) $
converge in the $C^{\infty }$ topology for $t>0$ we see that $g\left(
t\right) $ has positive sectional curvature for $t>0.$ This shows that the
metric is irreducible. Below we will show that the product metric on $%
M\times \mathbb{R}^{2}$ has nonnegative isotropic curvature for $t>0.$ The
Brendle-Schoen classification (see \cite[Theorem 2]{Brendle-SchoenII}) of
such metrics then shows that $M$ and hence $M_{i}$ (for sufficiently large $%
i $) are diffeomorphic to a sphere or compact rank one symmetric space,
giving the required contradiction.

In Brendle-Schoen \cite{Brendle-SchoenI} the authors show that a
quarter-pinched metric also has isotropic curvature even after the metric
has been multiplied with the factor $\mathbb{R}^{2}.$ The proof is entirely
algebraic and also shows that if the metric is almost quarter-pinched, then
there is a small lower bound for the isotropic curvatures, again even after
we have added the $\mathbb{R}^{2}$ factor. Thus the metrics $\left(
M_{i},g_{i}\left( 0\right) \right) $ have a lower bound $-\varepsilon _{i}$
for the isotropic curvatures on $M_{i}\times \mathbb{R}^{2},$ where $%
\varepsilon _{i}\rightarrow 0$ as $i\rightarrow \infty .$ (The $\varepsilon
_{i}$ can in fact be taken to be a multiple of $\frac{1}{i},$ but this is
not important for the proof.)

Finally we need to estimate the lower bound for the isotropic curvatures of $%
g_{i}\left( t\right) $ on $M_{i}\times \mathbb{R}^{2}.$ It is a general fact
that they will be bounded from below by $-\varepsilon _{i}\exp \left(
Ct\right) ,$ where $C$ is a constant that depends only on the curvature
bounds for $g_{i}\left( 0\right) .$ This was established by Hamilton in the
proof of \cite[Theorem 4.3]{HamiltonIII} as a general property of solutions
to heat flows. Specifically consider a heat flow%
\begin{equation*}
\partial _{t}T=\Delta T+\phi \left( T\right)
\end{equation*}%
on tensors $T$ and a continuous, compact, and convex condition $X$ on the
tensors $T$ that is also invariant under the natural $O\left( n\right) $
action induced on tensors and future invariant under the ODE $\partial
_{t}T=\phi \left( T\right) $. Hamilton shows that $X$ is then also invariant
under the heat flow by showing that if a solution to the heat flow starts
out at some distant $\varepsilon >0$ from $X,$ then at time $t$ it'll be at
most distance $\varepsilon \exp \left( Ct\right) $ away from $X.$ Here $C$
depends on the set $X$ and $\phi .$ The condition $X$ is rarely compact in
applications, but this can easily be finessed by modifying the PDE and ODE,
as we shall now discuss. In our situation we are considering the evolution
of the curvature tensor%
\begin{equation*}
\partial _{t}R=\Delta R+R^{2}+R^{\#}
\end{equation*}%
under the Ricci flow. Initially the curvature tensor is bounded by some
fixed constant $K_{1}$ that depends only on pinching constants and
dimension. This shows that the flow exists on some fixed time interval $%
\left[ 0,t_{0}\right] $ and that the curvature tensors are bounded by some
other controlled constant $K_{2}$ on this interval. If we have a continuous,
convex, $O\left( n\right) ,$ and future ODE invariant condition $X,$ then we
can intersect it with the compact, convex, and $O\left( n\right) $ invariant
set: $\left\vert R\right\vert \leq K,$ where $K>K_{2}.$ Next we modify the
vector field $\phi \left( R\right) =R^{2}+R^{\#}$ outside the region $%
\left\vert R\right\vert \leq K_{2}$ to make the new condition invariant
under the ODE $\partial _{t}R=R^{2}+R^{\#}.$ As long as we only consider
initial conditions where the curvature is bounded by $K_{1}$ the solutions
will not be affected by these changes.

Finally we can now use that for $t>0$ the metrics $g_{i}\left( t\right) $
converge to $g\left( t\right) $ in the $C^{\infty }$ topology to see that
the limit metrics $g\left( t\right) $ have nonnegative isotropic curvatures
on $M\times \mathbb{R}^{2}$, as desired. The proof of Theorem \ref{main} is
now complete.


\begin{thebibliography}{99}
\bibitem{Abresch-MeyerI} U. Abresch and W.T. Meyer, \emph{Pinching below }$%
1/4$\emph{, injectivity radius, and conjugate radius}. J. Differential Geom.
40 (1994), no. 3, 643--691.

\bibitem{Abresh-MeyerII} U. Abresch and W.T. Meyer, \emph{A sphere theorem
with a pinching constant below }$1/4$. J. Differential Geom. 44 (1996), no.
2, 214--261.

\bibitem{Berger} M. Berger, \emph{Sur les vari\'{e}t\'{e}s riemanniennes pinc%
\'{e}es juste au-dessous de 1/4}. Ann. Inst. Fourier (Grenoble) 33 (1983),
no. 2, 135--150.

\bibitem{Bohm-Wilking} C. B\"{o}hm and B. Wilking. \newblock Manifolds with
positive curvature operators are space forms. \newblock To appear in Ann. of
Math. \newblock arXiv:math.DG/0606187.

\bibitem{Brendle-SchoenI} S. Brendle and R. Schoen, \emph{Manifolds with
1/4-pinched Curvature are Space Forms}, arXiv: math.DG/0705.076

\bibitem{Brendle-SchoenII} S. Brendle and R. Schoen, \emph{Classification of
manifolds with weakly 1/4-pinched curvatures}. arXiv: math.DG/0705.3963. To
appear in Acta Math.

\bibitem{Chow-Knopf} B. Chow and D. Knopf. \newblock The Ricci flow: an
introduction. \newblock{\em Mathematical Surveys and Monographs}, vol. 110,
AMS, Providence, RI, 2004.

\bibitem{Chow-Lu-Ni} B. Chow, P. Lu, and L. Ni. \newblock Hamilton's Ricci
flow. \newblock {\em Graduate studies in Mathematics}, AMS, Providence, RI,
2006.

\bibitem{Durumeric} O. Durumeric, \emph{A generalization of Berger's theorem
on almost 1/4-pinched manifolds}. II. J. Differential Geom. 26 (1987), no.
1, 101--139.

\bibitem{HamiltonII} R. Hamilton. \newblock Three-manifolds with positive
Ricci curvature. \newblock { \em J. Diff. Geom.} 17(1982), 255-306.

\bibitem{HamiltonIII} R. Hamilton. \newblock Four-manifolds with positive
curvature operator. \newblock {\em J. Diff. Geom.} 24(1986), no. 2, 153-179

\bibitem{HamiltonIV} R. Hamilton. \newblock The Formation of Singularities
in the Ricci Flow. 
\newblock {\em Surveys in Differential
Geometry, Vol 2. (Cambridge, MA, 1993)}, Int. Press, Cambridge, MA, 1995,
7-136.

\bibitem{Petersen} P. Petersen, \emph{Riemannian Geometry 2}$^{nd}$\emph{\ Ed%
}, New York, Springer Verlag, 2006.

\bibitem{Petersen-Wilhelm} P. Petersen and F. Wilhelm, \emph{An exotic
sphere with positive sectional curvature}, arXiv:math.DG/0805.0812

\bibitem{Rong} X. Rong, \emph{On the fundamental groups of manifolds of
positive sectional curvature}. Ann. of Math. (2) 143 (1996), no. 2, 397--411.
\end{thebibliography}
\end{document}